\documentclass[a4paper]{amsart}

\usepackage{amssymb,stmaryrd,amsmath,graphics,verbatim}
\usepackage{latexsym}
\usepackage{eucal}
\makeatletter
\@addtoreset{equation}{section}\makeatother

\def\Car{1\hskip-.9mm{\rm l}}

\def\N{\mathbb{N}}
\def\Re{\mathbb{R}}
\def\Ric{\mathop{\rm Ric}\nolimits}
\def\S{\mathbb{S}}

\def\Diam{\mathop{\rm Diam}\nolimits}
\def\Vol{\mathop{\rm Vol}\nolimits}

\def\o {\circ}

\newtheorem{theorem}{Theorem}[section]
\newtheorem{lemma}[theorem]{Lemma}

\newtheorem{corollary}[theorem]{Corollary}

\newtheorem{remark}[theorem]{Remark}

\begin{document}
\title[Simplicial Volume in almost non negative Ricci curvature]{Bounds on the volume entropy and simplicial volume in Ricci curvature $L^p$-bounded from below}
\author{Erwann Aubry}
\address{Laboratoire Dieudonn\'e Univ. Nice Sophia-Antipolis Parc Valrose 06108 Nice FRANCE}
\email{eaubry@math.unice.fr}

\begin{abstract}
Let $(M,g)$ be a compact manifold with Ricci curvature almost bounded from below and $\pi:\bar{M}\to M$ be a normal, Riemannian cover. We show that, for any non-negative function $f$ on $M$, the means of $f\o\pi$ on the geodesic balls of $\bar{M}$  are comparable to the mean of $f$ on $M$. Combined with logarithmic volume estimates, this implies bounds on several topological invariants (volume entropy, simplicial volume, first Betti number and presentations of the fundamental group) in Ricci curvature $L^p$-bounded from below.
\end{abstract}


\maketitle

\section{Introduction}

We say that a compact manifold has {\sl Ricci curvature $L^p$ bounded from below (denoted  R.Lp.b. subsequently) by k} when the following quantity is small
$$\displaylines{\hfill \|\rho_k\|_p=\Bigl(\frac{1}{\Vol M}\int_{M}\rho_k^p\Bigr)^\frac{1}{p},\hfill}$$
where $k$ is a real, $\rho_k=\bigl(\underline{\Ric}-k(n-1)\bigr)_-$, $\underline{\Ric}(x)$ is the least eigenvalue of the Ricci tensor $\Ric$ at $x$ and $f_-=\max(-f,0)$.
Our purpose is to study some topological constraints under such control of the curvature.

On manifolds whose Ricci curvature is bounded from below in the usual sense, topological bounds are often derived from geometric comparison theorems on the volume of geodesic balls or spheres applied to some Riemannian covers (see for instance the proofs of the bounds on the first Betti number $\beta_1$ or the simplicial volume due to M. Gromov \cite{Gr,Gr2} or the bounds on the fundamental group due to M. Anderson \cite{An}).

 There are many extensions of the classical volume comparison theorems to manifolds with {\sl R.Lp.b.} (see for instance \cite{Ga2,Y,PSW,PW1,PS,Au6}), unfortunately, since these $L^p$ lower bounds on the Ricci curvature were not known to be preserved under non-finite Riemannian coverings, these volume estimates led to somewhat unsatisfactory topological bounds. For instance, in \cite{Ga2,RY} bounds on the simplicial volume or volume entropy in {\sl R.Lp.b.} are derived under additional assumptions on the universal Riemannian cover (such as a lower bound on the systole or a $L^p$ lower bound on the Ricci curvature of the cover itself) or on the curvature (such as a $L^p$ bound on the sectional curvature or an $L^\infty$ bound on the Ricci curvature).

To avoid this problem, S. Gallot developed  in \cite{Ga2} an approach based on M\"oser iterations which allows to bound the  harmonic topological invariants (i.e. those which, as the Betti numbers or the $\hat{A}$-genus, are associated to the multiplicities of certain eigenvalues of some natural Schr\"odinger operators). Gallot's bounds depend only on the control of some Sobolev constants of the manifold, and so does not rely on a control of the Ricci curvature on some covering space. The Gallot's bound on the first Betti number were, until recently, the only topological bound to rely exclusively on a $L^p$ lower bound of the Ricci curvature (this approach was also used in \cite{RY} for instance).

In \cite{Au6}, by constructing Dirichlet domains with multiplicity and using volume comparison results for star-shaped domains, we were able to get a $L^p$-lower bound on the Ricci curvature of the universal cover of any manifold with almost positive Ricci curvature. This was designed to extend the Myers's finiteness of the fundamental group in almost positive Ricci curvature (this was the first, non-harmonic bound on the topology of a manifold with {\sl R.Lp.b.}). Since, Z. Hu and S. Xu \cite{HX} have used the same technique to get some extensions of the Anderson's bound on the presentation of the $\pi_1$ in {\sl R.Lp.b.}  They also show that the technique developed in \cite{Au6} combined with the schema of proof developed by M.Gromov  \cite{Gr} in $\Ric\geq k$ allow to recover the Gallot's bound \cite{Ga2} on the first Betti number in {\sl R.Lp.b.}.

Generalizing the arguments used in \cite{Au6}, we prove the following lemma.

\begin{lemma}\label{Fundlemma}
 Let $n\geq2$ be an integer and $p>n/2$ be a real.  There exists a constant $\zeta(p,n,D\sqrt{|k|})>0$ such that if $(M^n,g)$ satisfies $\Diam M\leq D$ and $D^2\|\rho_k\|_p\leq\zeta(p,n,\sqrt{|k|}D)$ (for some $k\leq0$) then for any non-negative function $f$ on $M$, for any normal Riemannian cover $\pi:(\overline{M},\bar{g})\to (M,g)$, for any $\bar{x}\in\overline{M}$ and for any $R\geq 3D$, we have that
$$\displaylines{ \frac{1}{3^{n+1}e^{2(n-1)\sqrt{|k|}D}}\frac{1}{\Vol M}\int_{M}f\leq\frac{1}{\Vol B_{\bar{x}}(R)}\int_{B_{\bar{x}}(R)}f\o\pi\leq3^{n+1}e^{2(n-1)\sqrt{|k|}D}\frac{1}{\Vol M}\int_{M}f.}$$
\end{lemma}

 Applied to $f=\rho_k$, this lemma implies that {\em any normal Riemannian cover of a compact manifold with R.Lp.b. is also with R.Lp.b.} (where in the case of a non compact cover space $\overline{M}$, we say that $(\overline{M},\bar{g})$ is {\sl R.Lp.b.} if for a given $R$, $\sup_{\bar{x}\in \overline{M}}\int_{B_{\bar{x}}(R)}\rho_k^p/\Vol B_{\bar{x}}(R)$ is small).
Lemma \ref{Fundlemma} combined with some comparison results on logarithmic derivative of the volume of balls in {\sl R.Lp.b.} (see Theorem \ref{m3}), allows to extend the Gromov's bounds on the volume entropy and simplicial volume to manifold with Ricci curvature $L^p$-bounded from below.

\subsection{Bound on the volume entropy}

For any manifold $(M^n,g)$, we denote by $(\widetilde{M},\tilde{g})$ the Riemannian universal cover. The volume entropy of $(M^n,g)$ is defined by
$${\rm Ent}(M)=\limsup_{R\to+\infty}\frac{\ln\bigl( \Vol B_{\tilde{x}}(R)\bigr)}{R},$$
where $\tilde{x}$ is any point of $\tilde{M}$.
We set $A_k(R)$ (resp. $L_k(R)$) the volume of a geodesic ball (resp. a geodesic sphere) of radius $R$ in the space form of constant curvature $k$.

\begin{theorem}\label{Volsimp}
Let $n\geq2$ be an integer and $p>n/2$, $k\leq0$ be some reals. There exist some constants $\zeta(p,n,\sqrt{|k|}D)>0$ and $C(p,n)>0$ such that if $(M^n,g)$ satisfies $\Diam(M)\leq D$ and $D^2\|\rho_k\|_p\leq \zeta(p,n,\sqrt{|k|}D)$ then
$${\rm Ent}(M)\leq C(p,n)e^{\frac{(n-1)\sqrt{|k|}D}{p}}\Bigl(\frac{1}{\Vol M}\int_M\underline{\Ric}_-^p\Bigr)^\frac{1}{2p}.$$
\end{theorem}

\begin{remark}
  The constant $C(p,n)$ is equal to 
$$2^\frac{1}{2p}(n-1)^\frac{p-1}{2p}\Bigl(\frac{4p(p-1)}{(2p-1)(2p-n)}\Bigr)^\frac{p-1}{2p}3^\frac{n+1}{2p},$$
hence it tends to $\sqrt{n-1}$ when $p$ tends to $+\infty$, and so we recover the usual bound ${\rm Ent}(M)\leq(n-1)\sqrt{|k|}={\rm Ent}(\mathbb{H}^n,\frac{can}{k})$ when $\Ric\geq k(n-1)$.
\end{remark}

\begin{remark}
  In Theorem \ref{Volsimp}, we only assume a $L^p$ lower bound on the Ricci curvature of $M$. The bounds on the volume entropy in \cite{Ga2,RY} assume supplementary restrictions on the action of the fundamental group of $M$ on $\tilde{M}$.
\end{remark}

\subsection{Bounds on the simplicial volume}

Let $c$ a closed $l$-chain of $M$ and $\|[c]\|=\inf\{\sum_i|\alpha_i|/\,[\sum_i\alpha_ic_i]=[c]\}$, where the $c_i$ are elementary simplices. When $c$ is the fundamental class of $M$ then $\|[c]\|$ is called the simplicial volume of $M$ and denoted $\Vol_s(M)$. The bounds on the simplicial volume or norms of homology-classes by the the volume entropy due to M. Gromov \cite{Gr2} give the following corollaries.

\begin{corollary}\label{Vols1}
   Let $n\geq2$ be an integer and $p>n/2$, $k\leq0$ be some reals. There exist some constant $\zeta(p,n,D\sqrt{|k|})>0$ and $C(p,n)>0$ such that if $(M^n,g)$ satisfies $\Diam(M)\leq D$ and $D^2\|\rho_k\|_p\leq \zeta(p,n,\sqrt{|k|}D)$ then
$$\Vol_s(M)\leq C(p,n)^ne^\frac{n(n-1)\sqrt{|k|}D}{p}\Bigl(\frac{1}{\Vol M}\int_M\underline{\Ric}_-^p\Bigr)^\frac{n}{2p}\Vol(M),$$
and for any closed $l$-form we have that
$$\bigl\|[c]\bigr\|\leq l!C(p,n)^le^\frac{l(n-1)\sqrt{|k|}D}{p}\Bigl(\frac{1}{\Vol M}\int_M\underline{\Ric}_-^p\Bigr)^\frac{l}{2p}\Vol_l(c).$$
Once again, with the constant $C(p,n)$ computed in this paper we recover the Gromov's bounds when $p$ tends to $+\infty$.
\end{corollary}

We can also prove the following bound on the simplicial volume. It relies on other volume estimates and is not a corollary of the previous result

\begin{theorem}\label{Vols2}
   Let $n\geq2$ be an integer and $p>n/2$, $k\leq0$ be some reals. For any $\epsilon>0$, there exist a constant $\zeta(p,n,D\sqrt{|k|},\epsilon)>0$ such that if $(M^n,g)$ satisfies $\Diam(M)\leq D$ and $D^2\|\rho_k\|_p\leq \zeta(p,n,\sqrt{|k|}D,\epsilon)$ then
$$\Vol_s(M)\leq n!\Bigl(\frac{\Gamma(\frac{n}{2})}{\sqrt{\pi}\Gamma(\frac{n+1}{2})}\Bigr)^n\bigl((n-1)^n|k|^n+\epsilon\bigr)\Vol M,$$
and for any closed $l$-form we have that
$$\bigl\|[c]\bigr\|\leq l!\bigl((n-1)^l|k|^\frac{l}{2}+\epsilon\bigr)\Vol_l(c).$$
\end{theorem}

\subsection{Other bounds on the topology}

When $(M^n,g)$ is a manifold with {\sl R.Lp.b}, then Lemma \ref{Fundlemma} gives some $L^p$ lower bounds on the Ricci curvature of the Riemannian, normal covers of $M$. Using the various generalizations of the Bishop-Gromov inequality in {\sl R.Lp.b.} we can then generalize many topological bounds known when $\Ric\geq k(n-1)$ (see \cite{Au6,HX}). For instance, we can easily prove the following generalization of a result due to M.~Anderson.

\begin{theorem}\label{And2}
  Let $n\geq2$ be an integer and $p>n/2$, $k\leq0$, and $D,V>0$ be some reals. There exist some constants $\zeta(p,n,D\sqrt{|k|})>0$ and $N(p,n,D\sqrt{|k|})$ such that if $(M^n,g)$ satisfies $\Diam(M)\leq D$, $\Vol(M)\geq V$ and $\int_M\rho_k^p\leq V\cdot\zeta(p,n,D\sqrt{|k|})$ then any subgroup $\Gamma\subset\pi_1(M)$ generated by elements of length less than $\frac{VD}{N(p,n,D\sqrt{|k|})}$ has order less than $N(p,n,D\sqrt{|k|})/V$.
\end{theorem}

However, in all the extensions of the Bishop-Gromov theorem in {\sl R.Lp.b.} quoted above, to bound from above the quotient $\Vol B_{\tilde{x}}(R)/Vol B_{\tilde{x}}(r)$ (for $r\leq R$) we need to have
$$\frac{1}{\Vol B_{\bar{x}}(R)}\int_{B_{\bar{x}}(R)}(\underline{\Ric}-k)_-^p\leq C(p,n)R^{-2p}.$$
So, if extensions in {\sl R.Lp.b} of topological bounds such as the Gromov's bound on the first Betti number or the Anderson's bounds on the fundamental group (where the Bishop-Gromov estimate is only needed for balls of Radius comparable to the diameter of $M$) more or less readily follow from Lemma \ref{Fundlemma} and the already known estimates on the volume in {\sl R.Lp.b.}, on the contrary, when we need a Bishop-Gromov estimate for balls of arbitrary large radius in the cover (as in the proof of Theorem \ref{Volsimp} or in the proof of the Milnor's polynomial growth of the fundamental group in non negative Ricci curvature) the volume estimates of \cite{PSW,PW1,PS,Au6} are of no help. Note that to avoid this problem in their extension of the Milnor's polynomial growth in \cite{RY,HX}, the authors assume a supplementary lower bound on the systole or the volume of the manifold. This is a strong assumption since then it only remains a finite number of possible $\pi_1(M)$, which gives a bound on the radius of the balls we have to consider in the universal cover. In this paper, we prove a logarithmic estimate on the volume of balls (see theorem \ref{m2}) that applies for balls of arbitrary large radii (unfortunately, it is not good enough to imply a polynomial growth of $\Vol B_{\bar{x}}(R)$ in almost positive Ricci curvature).
\bigskip

In the last section of the paper, we collect some precompactness results for the Gromov-Hausdorff distance on the set of normal, Riemannian covers of compact manifolds with Ricci curvature $L^p$-bounded from below.
\medskip

\subsection{Remarks about our assumptions}

In \cite{Ga2}, S.Gallot constructs (for any $p>n/2$) a sequence of manifolds $(M_k,g_k)$ (example A.3 of the appendix) which satisfies
$$\displaylines{\hfill\Diam(g_k)\leq D,\hfill \frac{D^{2p}}{\Vol M}\int_M\bigl(\underline{\Ric}(g_k)\bigr)_-^p\leq C(p,n),\hfill b_1(M_k)\to+\infty.\hfill}$$
and can be easily adapted to also satisfy $\Vol_s(M_k)\to+\infty$ (since we can perform a connected sum of each $(M_k,g_k)$ with any fixed manifold $M$ without loosing the bounds on diameters and curvatures for $k$ large enough). Hence, if $D^2\|\rho_k\|_p$ is not assumed smaller than an universal constant, we can not bound the first Betti number,  the simplicial volume nor the volume entropy by the quantity $\frac{1}{\Vol M}\int_M\underline{\Ric}_-^p$.

Similarly, Examples A.2. of \cite{Ga2} shows that Theorem \ref{Volsimp} is not valid for $p\leq n/2$ (even if we replace $\frac{\Diam(M)^n}{\Vol M}\int_M\underline{\Ric}_-^{n/2}$ by $\int_M(\underline{\Ric})_-^{n/2}$). 

We proved in \cite{Au6} (Prop 9.3) that any compact manifold admit a metric with $D^2\|\rho_{1}\|_\frac{n}{2}$ as small as we want and so Corollary 1.5 and Theorem 1.6 are false in the case $p\leq\frac{n}{2}$.

Eventually, note that a still open conjecture due to M.Gromov asserts that $\Vol_s(M)$ should be bounded by a function of $\int_M(\underline{\Ric})_-^{n/2}$.
\medskip

\textbf{Acknowledgement} A part of these results were announced during the colloquium in honor of Marcel Berger {\sl Vari\'et\'es d'Einstein aujourd'hui et demain} (CIRM Nov. 2007). We thank S. Gallot, G.Wei and S.Rosenberg for their valuable comments on this work.

\section{Proof of Lemma \ref{Fundlemma}}\label{fundlemm}

The proof  of Lemma \ref{Fundlemma} relies on geometric group action arguments and some volume estimates in Ricci curvature $L^p$ bounded from below. We first recall some notations and results of \cite{Au6}.
\medskip

\noindent{\bf Notations.}  We denote by $U_x$ the {\em injectivity domain} at $x\in M$ and we identify points of $U_{x}{\setminus}\{0_{x}\}$ with their polar coordinates $(r,v)$ in $\Re^*_+{\times}\S_x^{n-1}$ (where $\S_x^{n-1}$ is the set of normal vectors at $x$). {\em We note $dv_g$ the Riemannian measure and set ${\rm exp}_{x}^*v_g=\theta(r,v)\,dr\, dv$, where $dv$ and $dr$ are the canonical measures of $\S_x^{n-1}$ and $\Re^*_+$}, and {\em we extend $\theta$ by $0$ to $\Re^*_+{\times}\S^{n-1}$}. {\em We denote by $h(r,v)$ the mean curvature at ${\rm exp}_x(rv)$ of the sphere centered at $x$ and of radius $r$ and we set} $\psi_k(r,v)=\bigl(\frac{(n-1)\sqrt{k}\cosh(\sqrt{k}r)}{\sinh(\sqrt{k}r)}{-}h(r,v)\bigr)_-$. 
Given $T$, a subset of $M$ star-shaped at $x$, let $A_T(r)$ be the volume of $B_x(r)\cap T$ and $L_T(r)$ be the $n-1$-volume of $\bigl(r\S^{n-1})\cap U_x\cap T_x$. We have $L_T(r)=\int_{\S^{n{-}1}_x}\Car_{T_x}\theta(r,v)dv$ and $A_T(r)=\int_0^rL_T(t)\,dt$. We set $L_k$, $A_k$ the corresponding functions on the space form of sectional curvature $k$ and $s_k(t)=\sinh(\sqrt{|k|}t)/\sqrt{|k|}$ if $k<0$, $s_0(t)=t$, $s_k(t)=\sin(\sqrt{k}t)/\sqrt{k}$ if $k>0$.

\begin{lemma}\label{accfinisL} $L_T$ is a right continuous, left lower semi-continuous function. $A_T$ is a continuous, right differentiable function of derivative $L_T$. Moreover, the function
$$f(r)=\frac{L_T(r)}{s_k(r)^{n-1}}-\int_0^r\int_{\S^{n-1}_x}\frac{\Car_{T_{x}}\psi_k\theta}{s_k(s)^{n-1}}\,dvds$$
is decreasing on $\Re^*_+$.
\end{lemma}

Moreover, the mean curvature can be controlled by integral of the Ricci curvature thanks to the following Lemma (see Lemma 4.1 in \cite{Au6}).

\begin{lemma}\label{fondcomp}
Let $p>n/2$, $r>0$ be some reals. We have
$$\displaylines{\psi_k^{2p-1}(r,v)\,\theta(r,v)\leq(2p{-}1)^p\left(\frac{n{-}1}{2p{-}n}\right)^{p{-}1}\int_0^{r}\rho_k^p(t,v)\theta(t,v)\,dt}$$
 for all normal vector $v\in\S^{n-1}_x$. Note that this inequality holds with $\theta$ replaced by $\Car_{[0,s[}\,\theta$.
\end{lemma}

\begin{remark}
Lemma \ref{fondcomp} holds in dimension $2$ with $p=1$ and the constant before $\int_0^r\rho_k^p\theta$ equal to $1$.
\end{remark}

\paragraph{Proof of Lemma \ref{Fundlemma}}

We set $s_k(t)=\frac{\sinh(\sqrt{k}t)}{\sqrt{k}}$ for $k<0$ and $s_0(t)=t$.
By geometric action group arguments, we prove in \cite{Au6} that there exists a subset  $T\subset\overline{M}$ such that
\begin{itemize}
\item $B_{\bar{x}}(R)\subset T\subset B_{\bar{x}}(R+D)$,
\item $T$ is star shaped at $\bar{x}$,
\item if $\Gamma$ is the deck transformation group of $\pi:\overline{M}\to M$, then $x\mapsto \#\bigl(T\cap\pi^{-1}(x)\bigr)$ is constant on $M$, and so we have for any function $f$ on $M$
$$\frac{1}{\Vol T}\int_T(f\o\pi)\,dv_{\bar{g}}=\frac{1}{\Vol M}\int_{M}f\,dv_g,$$
\end{itemize}
By Lemma \ref{fondcomp}, we have for any $R-D\leq t\leq r\leq R+D$, that
$$\displaylines{s_k^{n-1}(t)L_T(r)\leq s_k^{n-1}(r)L_T(t)+s_k^{n-1}(r)\int_t^r\int_{\S^{n-1}_{\bar{x}}}\Car_{T_{\bar{x}}}\bar{\psi}_k\bar{\theta}\,dvds\hfill\cr
\leq s_k^{n-1}(r)L_T(t) +s_k^{n-1}(r)(A_T(r)-A_T(t))^\frac{2(p-1)}{2p-1}(\int_t^r\int_{\S^{n-1}_{\bar{x}}}\Car_{T_{\bar{x}}}\bar{\psi}_k^{2p-1}\bar{\theta}\,dvds)^\frac{1}{2p-1}\cr
\leq s_k^{n-1}(r)\Bigl(L_T(t)+C_1(p,n)A_T(R+D)D^\frac{1}{2p-1}(\int_{T}\frac{\bar{\rho}_k^p}{\Vol T})^\frac{1}{2p-1}\Bigr),}$$
where we have also used the H\"older inequality and Lemma \ref{fondcomp}, $\Vol T\leq A_T(R+D)$ and where $C_1(p,n)$ is equal to $\bigl(\frac{(2p-1)^p(n-1)^{p-1}}{(2p-n)^{p-1}}\bigr)^\frac{1}{2p-1}$ and tends to $\sqrt{n-1}$ when $p$ tends to $+\infty$.
 By integration of this inequality between $R- D$ and $R$  with respect to $t$ and between $R$ and $R+D$ with respect to $r$, we get that
$$\displaylines{\bigl(A_T(R+D)-A_T(R)\bigr)\int_{R-D}^Rs_k^{n-1}\hfill\cr
\hfill\leq\Bigl(A_T(R)+C_1(p,n)A_T(R+D)\bigl(D^2\|\rho_k\|_p\bigr)^\frac{p}{2p-1}\Bigr)\int_R^{R+D} s_k^{n-1},}$$
which implies that
$$\Bigl(\int_{R-D}^R s_k^{n-1}-C_1(p,n)\bigl(D^2\|\rho_k\|_p\bigr)^\frac{p}{2p-1}\int_R^{R{+}D} s_k^{n-1}\Bigr)\frac{A_T(R{+}D)}{A_T(R)}\leq\int_{R-D}^{R+D} s_k^{n-1}.$$
By convexity of the function $A_k$, for any $R\in[3D,+\infty[$, we have
$$\frac{\int^{R+D}_{R-D}s_k^{n-1}}{\int^R_{R-D}s_k^{n-1}}\leq\frac{2D}{D}\Bigl(\frac{s_k(R+D)}{s_k(R-D)}\Bigr)^{n-1}\leq2\Bigl(\frac{s_k(3D)}{s_k(D)}\Bigr)^{n-1}\leq 3^ne^{2(n-1)\sqrt{|k|}D}$$
So, if $D^2\|\rho_k\|_p\leq\Bigl[\frac{1}{3^{n+1}C_1e^{2(n-1)\sqrt{|k|}D}}\Bigr]^\frac{2p-1}{p}$ then $\frac{A_T(R+D)}{A_T(R)}\leq3^{n+1}e^{2(n-1)\sqrt{|k|}D}$. Since we have $A_T(R+D)=\Vol T$ and $A_T(R)=\Vol B_{\bar{x}}(R)$, we get
$$\displaylines{\frac{1}{\Vol B_{\bar{x}}(R)}\int_{B_{\bar{x}}(R)}\!\!(f\o\pi)\leq\frac{A_T(R+D)}{A_T(R)}\int_T\frac{(f\o\pi)}{\Vol T}\leq3^{n+1}e^{2(n-1)\sqrt{|k|}D}\int_M\frac{f}{\Vol M}.}$$
To prove the inverse estimate, consider $T'$ the star shaped subset constructed as above for the set $B_{\bar{x}}(R-D)$. Then we have $T'\subset B_{\bar{x}}(R)$ and
$$\frac{1}{\Vol B_{\bar{x}}(R)}\int_{B_{\bar{x}}(R)}(f\o \pi)\geq\frac{\int_{T'}(f\o\pi)}{\Vol T'}\frac{\Vol B_{\bar{x}}(R-D)}{\Vol B_{\bar{x}}(R)}$$
and the quotient $\frac{\Vol B_{\bar{x}}(R-D)}{\Vol B_{\bar{x}}(R)}$ can be bounded from below as above.

\section{bound on the fundamental group}

Combining Lemma \ref{Fundlemma} with the Bishop-Gromov inequalities of  \cite{PW1} or \cite{Au6}, we get the following corollary.

\begin{corollary}\label{Bish-Gro2}
  Let $n\geq2$ be an integer and $p>n/2$, $k\leq0$, $\alpha>1$ be some reals. There exist some constants $\zeta(\alpha D\sqrt{|k|},p,n)>0$ and $C(\alpha D\sqrt{|k|},p,n)>0$ such that if $(M^n,g)$ satisfies $\Diam(M)\leq D$ and $D^2\|\rho_k\|_p\leq\zeta(D\sqrt{|k|},p,n)$ then for any normal cover $\overline{M}\to M$, for any $\bar{x}\in\overline{M}$, for any $\alpha D\geq R\geq r>0$, we have that
$$\displaylines{\frac{\Vol B_{\bar{x}}(R)}{\Vol B_{\bar{x}}(r)}\leq\frac{A_{k}(R)}{A_{k}(r)}\bigl(1+C(\alpha D\sqrt{|k|},p,n)(D^2\|\rho_{k}\|_{p})^\frac{p}{10}\bigr),\cr
\Vol B_{\bar{x}}(R)\leq A_{k}(R)\bigl(1+C(\alpha D\sqrt{|k|},p,n)(D^2\|\rho_{k}\|_{p})^\frac{p}{10}\bigr).}$$
\end{corollary}

This last corollary implies Theorem \ref{And2}. Indeed, we have $D^{2p}\|\rho_k\|_p^p\leq\zeta(p,n,D\sqrt{|k|})$, hence for $\zeta$ small enough, Corollary \ref{Bish-Gro2} applies for $\alpha=2$. We set
$$\displaylines{N(p,n, D\sqrt{|k|})=2A_{k}(2D)\bigl(1+C(2D\sqrt{|k|},p,n)(D^2\|\rho_{k}\|_{p})^\frac{p}{10}\bigr)\cr
L=\frac{DV}{N}.}$$
Let $\Gamma$ be a subgroup of $\pi_1(M)$ generated by $k$ elements $(g_1,\cdots,g_k)$ of lengths less than $L$. Let $(\tilde{M},\tilde{g})$ be the Riemannian, universal cover of $(M,g)$, let $\tilde{x}$ be a point of $\tilde{M}$ and $F$ be the Dirichlet domain of $\tilde{x}$. Then $\Diam (F)\leq D$ and if we set $\Gamma(N)=\{g\in\Gamma/|\gamma|\leq\frac{N}{V}\}$ then by Corollary \ref{Bish-Gro2} (with $r=0$), we have
$$\displaylines{\#\Gamma(N)\Vol F=\Vol\bigl(\bigcup_{\gamma\in\Gamma(N)}\gamma F\bigr)\leq\Vol\bigl(B_{\tilde{x}}(2D)\bigr)\leq \frac{N}{2}<\frac{N}{V}\Vol F}$$
hence $\# \Gamma(N)<N$. We infer that $\Gamma(N)=\Gamma$ and so $\#\Gamma\leq\frac{N}{V}$.

\section{Bound on the volume entropy}
It follows easily from the definition of the volume entropy that
$${\rm Ent}(M)\leq\sup_{\tilde{x}\in\widetilde{M}}\limsup_{R\to+\infty}\frac{\Vol S_{\tilde{x}}(R)}{\Vol B_{\tilde{x}}(R)}$$

To prove Theorem \ref{Volsimp}, we need an asymptotic logarithmic upper bound of the volume of balls.  We will use the following estimate due to S. Gallot (\cite{Ga2} Theorem 1, p.195 with $s=0$ and $H=S_x(R_0)$).

\begin{theorem}[\cite{Ga2}]\label{Ga}
  Let $n\geq 2$  be an integer and $p>n/2$, $k\leq0$ be some reals. There exist a constant $C(p,n)>0$ such that on any manifold $(M^n,g)$ with non finite volume, we have
$$\displaylines{\sup_{x\in M}\limsup_{R\to\infty}\frac{\Vol S_{x}(R)}{\Vol B_{x}(R)}\leq C(p,n)\limsup_{R\to \infty}\Bigl(\frac{1}{\Vol B_x(R)}\int_{B_x(R)}(\underline{\Ric}_-)^p\Bigr)^\frac{1}{2p},}$$
where $C(p,n)=2^\frac{1}{2p}(n-1)^\frac{p-1}{2p}\Bigl(\frac{4p(p-1)}{(2p-1)(2p-n)}\Bigr)^\frac{p-1}{2p}$ tends to $\sqrt{n-1}$ when $p$ tends to $+\infty$.
\end{theorem}

To prove Theorem \ref{Volsimp}, we just apply Theorem \ref{Ga} to the universal cover (which can be suppose with non finite volume) and bound the right-hand side by Lemma \ref{Fundlemma}.

\section{Bounds on the simplicial volume}\label{Hvol}

The bounds on the co-homology norms and simplicial volume of Corollary \ref{Vols1} follow from the following result.
\begin{theorem}[Gromov\cite{Gr}]
$$\displaylines{\|[c]\|\leq l!\bigl(\alpha(M)\bigr)^l\Vol_l(c),\cr
\Vol_{s}(M)=\|[M]\|\leq\bigl(\frac{\Gamma(\frac{n}{2})}{\sqrt{\pi}\Gamma(\frac{n+1}{2})}\bigr)^nn!\bigl(\alpha(M)\bigr)^n\Vol(M),}$$
where $\alpha(M)=\sup_{\tilde{x}\in\tilde{M}}\inf_{R>0}\frac{\Vol S_{\tilde{x}}(R) }{\Vol B_{\tilde{x}}(R)}$.
\end{theorem}

Since $\alpha(M)\leq {\rm Ent}(M)$, Theorem \ref{Vols1} follows from Theorem \ref{Volsimp}. The bound of Theorem \ref{Vols1} is optimal for large $p$ since it implies the optimal Gromov's bound when $p$ tends to $+\infty$. Theorem \ref{Vols2} is optimal when $\|\rho_k\|_p$ tends to $0$ for the same reason. To prove Theorem \ref{Vols2} we need another upper bound on the quotient $\Vol S_x(R)/\Vol B_x(R)$ that is optimal for small $\|\rho_k\|_p$.

\begin{theorem}\label{m2}
   Let $n\geq2$ be an integer and $p>n/2$, $k\leq0$ be some reals. There exists a constant $C(p,n)>0$ such that for any $(M^n,g)$, any $x\in M$ and any $R>0$, we have
$$\displaylines{\frac{\Vol S_x(R)}{\Vol B_x(R)}\leq\frac{L_{k}(R)}{A_{k}(R)}\Bigl(1+C(p,n)\bigl(R^{2p}\frac{\int_{B_x(R)}\rho_{k}^p}{\Vol B_x(R)}\bigr)^\frac{1}{2p-1}\Bigr).}$$
\end{theorem}

\begin{remark}
  Note that contrarily to other volume estimates in {\sl R.Lp.b}, the previous one does not requires smallness of $\|\rho_k\|_p$ to apply. It will have interesting consequences in the precompactness results of the last section.
\end{remark}

\begin{proof}
We set $\|\rho_k\|_{p,r}=\frac{\int_{B_x(r)}\rho_{k}^p}{\Vol B_x(r)}$. As in the proof of Lemma \ref{Fundlemma}, for any $t\leq r$ we have that
$$\displaylines{s_{k}^{n-1}(t)L_T(r)
\leq s_{k}^{n-1}(r)L_T(t)+C(p,n)^\frac{1}{2p-1}s_{k}^{n-1}(r)A_T(r)(r-t)^\frac{1}{2p-1}\|\rho_{k}\|_{p,r}^\frac{p}{2p-1},}$$
where $C(p,n)=(2p{-}1)^p\left(\frac{n{-}1}{2p{-}n}\right)^{p{-}1}$.
Integrating with respect to $t$ between $0$ and $r$, it gives us that
$$\frac{L_T(r)}{A_T(r)}\leq\frac{L_k(r)}{A_k(r)}\Bigl(1+C(p,n)^\frac{1}{2p-1}\frac{2p-1}{2p}\bigl(r^2\|\rho_k\|_{p,r}\bigr)^\frac{p}{2p-1}\Bigr).$$
\end{proof}

Combining Lemma \ref{Fundlemma} with the previous theorem (for large $R$) and the volume estimates on $L(R)$ of \cite{Au6} (for $R\leq 3D$), we get the following corollary, which implies Theorem \ref{Vols2}.

\begin{corollary}\label{m3} Let $n\geq2$ be an integer and $p>n/2$, $k\leq0$ be some reals. There exist some constants $\zeta(p,n,\sqrt{|k|}D)>0$ and $C(p,n)>0$ such that if $(M^n,g)$ satisfies $\Diam(M)\leq D$ and $D^2\|\rho_k\|_p\leq\zeta(p,n,\sqrt{|k|}D)$ then for any normal cover $\overline{M}\to M$, for any $\bar{x}\in\overline{M}$, for any $R>0$ and for any $k'\leq0$, we have that
$$\displaylines{\frac{\Vol S_{\bar{x}}(R)}{\Vol B_{\bar{x}}(R)}\leq\frac{L_{k'}(R)}{A_{k'}(R)}\Bigl(1+C(p,n)e^{\frac{2(n-1)\sqrt{|k|}D}{2p-1}}\bigl(\max(R^2,D^2)\|\rho_{k'}\|_{p}\bigr)^\frac{p}{2p-1}\Bigr).}$$
\end{corollary}

\section{Precompactness results}

Integrating the volume estimate of Theorem \ref{m2} we get the following Bishop-Gromov kind estimate.
\begin{corollary}\label{m4}
   Let $n\geq2$ be an integer and $p>n/2$, $k\leq0$ be some reals. There exists a constant $C(p,n)>0$ such that for any $(M^n,g)$, any $x\in M$ and any $R\geq r>0$, we have
$$\displaylines{\frac{\Vol B_x(r)}{\Vol B_x(R)}\leq\Bigl(\frac{A_{k}(r)}{A_{k}(R)}\Bigr)^{1+C(p,n)f(r,R)},}$$
where $f(r,R)=\sup_{t\in[r,R]}\bigl(t^{2p}\frac{\int_{B_x(t)}\rho_{k}^p}{\Vol B_x(t)}\bigr)^\frac{1}{2p-1}$.
\end{corollary}

As a by product of the previous volume estimates, we get the following precompactness result.
\begin{theorem}\label{Precomp}
Let $n\geq 2$ be an integer and $p>n/2$, $D>0$ be some reals. Let $f$ be any positive, locally bounded function on $]0,+\infty[$. The set of compact manifolds which satisfy
$$\displaylines{\hfill\sup_{x\in M}\frac{\int_{B_x(R)}\bigl(\underline{\Ric}-k(n-1)\bigr)_-^p}{\Vol B_x(R)}\leq f(R),\hfill and\hfill\Diam(M)\leq D\hfill}$$
for any $R>0$ is precompact for the Gromov-Hausdorff distance.
\end{theorem}
\begin{remark}
  Contrarily to the precompactness result of \cite{PW1} we do not have to suppose $f$ sufficiently small. Of course this result is more interesting when $f$ is chosen non bounded near $0$.
\end{remark}

 The same result holds for complete manifolds and the pointed Gromov-Hausdorff topology.
\begin{theorem}\label{Precompl}
Let $n\geq 2$ be an integer, $p>n/2$ be a real and $f$ be any positive, locally bounded function on $]0,+\infty[$. If $(M_i,x_i,g_i)_{i\in\N}$ is a sequence of pointed manifolds  which satisfy
$$\displaylines{\hfill\sup_{x\in B(x_i,R_i)}\frac{\int_{B_{x}(R)}\bigl(\underline{\Ric}-k(n-1)\bigr)_-^p}{\Vol B_{x}(R)}\leq f(R),\hfill}$$
for any $R_i\geq R>0$ (where $R_i\to+\infty$), then a sub-sequence converges to a length space in the pointed Gromov-Hausdorff topology.
\end{theorem}
\medskip

Integrating the volume estimate of Corollary \ref{m3} we get the following volume estimate.

\begin{corollary}\label{Bish-Gro}
 Let $n\geq2$ be an integer and $p>n/2$, $k\leq0$ be some reals. There exist some constants $\zeta(p,n,\sqrt{|k|}D)>0$ and $C(p,n)>0$ such that if $(M^n,g)$ satisfies $\Diam(M)\leq D$ and $D^2\|\rho_k\|_p\leq\zeta(p,n,\sqrt{|k|}D)$ then for any normal cover $\overline{M}\to M$, for any $\bar{x}\in\overline{M}$, for any $R\geq r>0$ and any $k'\leq0$, we have that
$$\displaylines{\frac{\Vol B_{\bar{x}}(R)}{\Vol B_{\bar{x}}(r)}\leq\Bigl(\frac{A_{k'}(R)}{A_{k'}(r)}\Bigr)^{1+C(p,n)e^{\frac{2(n-1)\sqrt{|k|}D}{2p-1}}\bigl(\max(R^2,D^2)\|\rho_{k'}\|_{p}\bigr)^\frac{p}{2p-1}}.}$$
\end{corollary}

This corollary readily implies the following precompactness result for the normal covers in {\sl R.Lp.b.}

\begin{theorem}\label{Precompcov}
Let $n\geq 2$ be an integer, $k$ and $p>n/2$ be some reals. There exists a constant $\zeta(p,n,k)>0$ such that if $(M_i,g_i)_{i\in\N}$ is a sequence of manifolds which satisfy $\frac{(\Diam M_i)^{2p}}{\Vol M_i}\int_M\rho_k^p\leq\zeta(p,n,k)$, then, for any normal, Riemannian cover $\pi_i:(\overline{M}_i,\bar{g}_i)\to (M_i,g_i)$ and any $\bar{x}_i\in\overline{M}_i$, the sequence $(\overline{M}_i,\bar{x}_i,\bar{g}_i)$ admits a  sub-sequence that converges in the Gromov-Hausdorff pointed topology.
\end{theorem}

\end{document}